\documentclass[twocolumn]{article}
\usepackage{amssymb,mathrsfs,amsmath,amsfonts,graphicx,pdfpages,csquotes,amsthm}

\newtheorem{theorem}{Theorem}[section]
\newtheorem{defn}{Definition}[section]

\newtheorem{cor}{Corollary}[section]

\newcommand{\maps}{\colon}

\newcommand{\op}{\mathrm{op}}

\newcommand{\namedcat}[1]{\mathsf{#1}}

\newcommand{\Set}{\namedcat{Set}}

\newcommand{\Rel}{\mathsf{Rel}}

\newcommand{\colim}{\mathrm{colim}}
\newcommand{\Alg}{\mathsf{Alg}}
\newcommand{\Coalg}{\mathsf{Coalg}}

\makeatletter
\newcommand*{\relrelbarsep}{.386ex}
\newcommand*{\relrelbar}{%
  \mathrel{%
    \mathpalette\@relrelbar\relrelbarsep
  }%
}
\newcommand*{\@relrelbar}[2]{%
  \raise#2\hbox to 0pt{$\m@th#1\relbar$\hss}%
  \lower#2\hbox{$\m@th#1\relbar$}%
}
\providecommand*{\rightrightarrowsfill@}{%
  \arrowfill@\relrelbar\relrelbar\rightrightarrows
}
\providecommand*{\leftleftarrowsfill@}{%
  \arrowfill@\leftleftarrows\relrelbar\relrelbar
}
\providecommand*{\xrightrightarrows}[2][]{%
  \ext@arrow 0359\rightrightarrowsfill@{#1}{#2}%
}
\providecommand*{\xleftleftarrows}[2][]{%
  \ext@arrow 3095\leftleftarrowsfill@{#1}{#2}%
}
\makeatother

\usepackage{subfiles}
\usepackage{stackengine}
\usepackage{mathtools}
\usepackage[utf8]{inputenc}
\usepackage{tikz}
\usepackage[a4paper,top=3cm,bottom=3cm,inner=3cm,outer=3cm]{geometry}

\usepackage{comment}

\usepackage{graphicx}
\usepackage{adjustbox}
\usepackage[all,2cell]{xy}\UseAllTwocells\SilentMatrices
\definecolor{darkgreen}{rgb}{0,0.45,0}
\definecolor{myurlcolor}{rgb}{0,.45,.2}
\definecolor{mylinkcolor}{rgb}{.1,.1,.7}
\usepackage[colorlinks,citecolor=darkgreen,urlcolor=myurlcolor,linkcolor=mylinkcolor,final,hyperindex,linktoc=page,pagebackref]{hyperref}

\usepackage[capitalize]{cleveref}
\crefname{equation}{}{}
\crefname{item}{}{}
\usepackage{enumerate}

\usepackage{tikz}
\usepackage{tikz-cd}
\usetikzlibrary{backgrounds,circuits,circuits.ee.IEC,shapes,fit,matrix}

\usepackage{verbatim}
\usetikzlibrary{arrows,shapes}

\tikzstyle{simple}=[-,line width=2.000]
\tikzstyle{arrow}=[-,postaction={decorate},decoration={markings,mark=at position .5 with {\arrow{>}}},line width=1.100]
\pgfdeclarelayer{edgelayer}
\pgfdeclarelayer{nodelayer}
\pgfsetlayers{edgelayer,nodelayer,main}

\tikzstyle{none}=[inner sep=0pt]

\definecolor{lblue}{rgb}{0,250,255}
\tikzstyle{species}=[circle,fill=yellow,draw=black,scale=1.15]
\tikzstyle{transition}=[rectangle,fill=lblue,draw=black,scale=1.15]
\tikzstyle{inarrow}=[->, >=stealth, shorten >=.03cm,line width=1.5]
\tikzstyle{empty}=[circle,fill=none, draw=none]
\tikzstyle{inputdot}=[circle,fill=black,draw=black, scale=.25]
\tikzstyle{inputarrow}=[->,draw=purple, shorten >=.05cm]
\tikzstyle{simple}=[-,draw=black,line width=1.000]
\tikzstyle{dot}=[circle,fill=black,draw=black, scale=.4]
\tikzstyle{inarrow}=[->, >=stealth, shorten >=.03cm,line width=1.5]

\definecolor{joecolor(x11)}{rgb}{0.0, 0.5, 0.5}

\definecolor{purple(x11)}{rgb}{0.8, 0, 0.8}

\begin{document}
\title{Beyond Initial Algebras and Final Coalgebras}
\author{Clemens Kupke, Jade Master, Ezra Schoen}

\maketitle

\begin{abstract}
  We provide a construction of the fixed points of functors which may not be inital algebras or final coalgebras. For an endofunctor $F$, this fixpoint construction may be expressed as a pair of adjoint functors between $F$-coalgebras and $F$-algebras. We prove a version of the limit colimit coincidence theorem for these generalized fixed points. 
\end{abstract}

\section{Middle Fixed Points}
The Knaster-Tarski theorem provides a construction of least and greatest fixed points for a monotone function $f: L \to L$ on a complete lattice $L$. Consider the following monotone function on the lattice $([0,1],\leq)$ of the interval of real numbers with the usual ordering. 
\begin{center}
 \includegraphics[scale=.3]{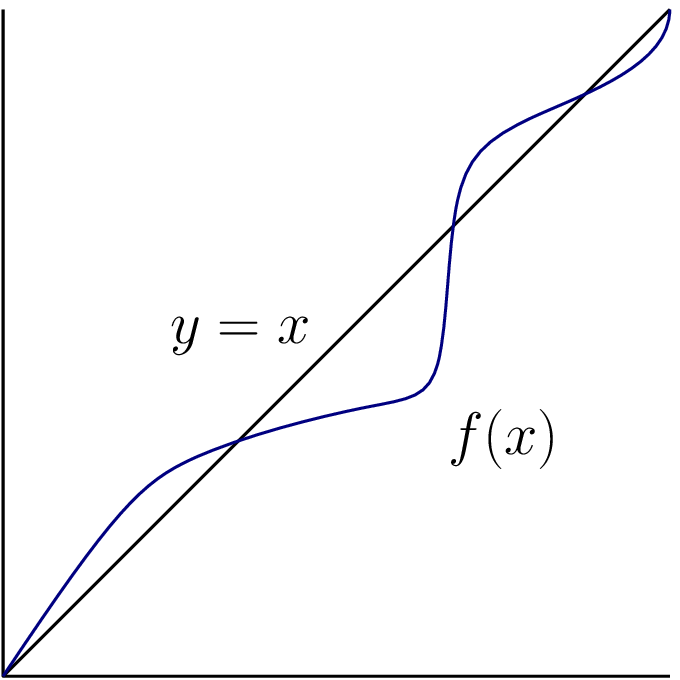}
 \end{center}
The function $f$ is overlayed with the function $y=x$. The intersection of the two curves indicate fixpoints of $f$. The least fixpoint of $f$ is $0$ and the greatest fixpoint is $1$ but there are $3$ other fixpoints in-between. These ``middle" fixpoints have a similar construction to the least and greatest ones. Given a ``pre-fixed point" i.e. a point $x \in [0,1]$ such that $x \leq f(x)$ we may find the first fixpoint above $x$ as
\[\mu(x) = \sup \{x,f(x),f^2(x),f^3(x),\ldots\} \]
where the $\ldots$ indicate iteration to a sufficiently large ordinal. Similarly, given a ``post-fixed point" $f(y) \leq y$, we may find the closest fixpoint below $y$
\[ \nu(y) = \inf \{ y ,f(y),f^2(y),f^3(y),\ldots\}.\] For a complete lattice $L$, let $Pre(f)$ be the suborder of $L$ consisting of only the pre-fixed points $x \leq f(x)$. Similarly, let $Post(f)$ be the suborder of post fixed points $f(y) \leq y$. Then there is a Galois connection
\[
\begin{tikzcd}
Pre(f) \ar[r,bend left,"\mu"] \ar[r,phantom,"\bot"] & Post(f) \ar[l,bend left,"\nu"]
\end{tikzcd}
\]
Being a Galois connection means that
\[\mu(x) \leq y \iff x \leq \nu(y) \]
In this abstract we will generalized this Galois connection to fixpoints of functors rather than monotone functions. When generalizing from posets to categories we make the replacements shown in Table \ref{table}. Although the `middle fixed point' construction is new, fixed points which are neither initial or final have been studied before. These are the rational fixed points of Ad{\'a}mek, Milius, and Velebil \cite{adamek2006iterative}.
\begin{figure*}[h]\label{table}
\begin{center}
\begin{tabular}{|c|c|}
\hline
    Poset & Category \\
    \hline
    Monotone Function $f$ & Functor $F$ \\
    Pre-fixed point of $f$ & $F$-coalgebra \\
    Post-fixed point of $f$ & $F$-algebra \\ $\sup \{f(x),f^2(x),f^3(x), \ldots \}$ & $\colim(X \to F(X) \to F^2(X) \to F^3(X) \ldots)$ \\
    $\inf\{f(x),f^2(x),f^3(x), \ldots \}$ &  $\lim (X \leftarrow F(X) \leftarrow F^2(X) \leftarrow F^3(X) \ldots )$ \\
    Galois connection & Adjunction\\
    \hline
\end{tabular}
\end{center}
\caption{Generalization of Posets to Categories}
\end{figure*}
\section{The Adjunction}
In what follows, $\omega$ will be the category finite ordinals i.e. your favorite category with a countable number of objects and one generating morphism from each object to the next. Pictorally, $\omega$ is the category
\[ \begin{tikzcd}\bullet \ar[r] & \bullet \ar[r] & \bullet \ar[r] &
\bullet \ldots \end{tikzcd}\]
\begin{theorem}\label{thm:1}
Suppose $C$ is a category with colimits of shape $\omega$ and limits of shape $\omega^{op}$ and suppose that $F: C \to C$ preserves limits and colimits of these shape. Then there is an adjunction
\[
\begin{tikzcd}
\Coalg(F) \ar[r,bend left,"\mu"] \ar[r,phantom,"\bot"]& \Alg(F) \ar[l,bend left,"\nu"]
\end{tikzcd}
\]
given by 
\[\mu (b \maps B \to FB) =\]\[\colim(\begin{tikzcd} B \ar[r,"b"] & F(B) \ar[r,"Fb"] & F^2(B) \ar[r,"F^2(b)"] & \cdots \end{tikzcd}) \]
\[\nu(a \maps FA \to A)=\]\[\lim(\begin{tikzcd} A & F(A) \ar[l,"a",swap] & F^2(A) \ar[l,"Fa",swap] & \ar[l,"F^2(a)",swap] \cdots \end{tikzcd} )\]
and defined on morphisms using the universal property of limits and colimits. 
\end{theorem}
\begin{proof}(Sketch) The adjunction isomorphism
\[\Alg(F)(\mu(b),a) \cong \Coalg(F)(b,\nu(a)) \]
for algebras $a \maps FA \to A$ and coalgebras $b \maps B \to FB$ relies on the fact that both sets are naturally isomorphic to the set of coalgebra to algebra homomorphisms from $b$ to $a$. Given a coalgebra to algebra homomorphism
\[
\begin{tikzcd}
    B \ar[d,"f"] \ar[r,"b"] & F(B)\ar[d,"F(f)"] \\
    A & \ar[l,"a"] F(A)
\end{tikzcd}
\]
we may iterate it countably many times to get a diagram
\[\begin{tikzcd}
	B & {F(B)} & {F^2(B)} & {F^3(B)\cdots} \\
	A & {F(A)} & {F^2(A)} & {F^3(A)\cdots}
	\arrow["b", from=1-1, to=1-2]
	\arrow["Fb", from=1-2, to=1-3]
	\arrow["{F^2b}", from=1-3, to=1-4]
	\arrow["a", from=2-2, to=2-1]
	\arrow["Fa", from=2-3, to=2-2]
	\arrow["{F^3a}", from=2-4, to=2-3]
	\arrow["{f}"', from=1-1, to=2-1]
	\arrow["F(f)", from=1-2, to=2-2]
	\arrow["F^2(f)", from=1-3, to=2-3]
	\arrow["F^3(f)\cdots", from=1-4, to=2-4]
\end{tikzcd}\]
The colimit of the top row is $\mu(b)$ and the maps going down and left form a cocone over the diagram for $\mu(b)$. Therefore the universal property for colimits induces a morphism $\mu(b) \to a$ and we state without proof that this is an algebra homomorphism. Similarly, the limit of the bottom row is the coalgebra $\nu(a)$ and the morphisms going right and down form a cone. The universal property of limits supplies a morphism $b \to \nu(a)$. We state without proof that this morphism is a coalgebra homomorphism and that the correspondences described here are natural in both arguments. 
\end{proof} 
\noindent Let $1$ be the terminal object of $C$ and let $0$ be the initial object. Then there is a unique algebra $1:F1 \to 1$ and $\nu(1)$ is the terminal coalgebra. Similarly, the initial algebra is given by $\mu(0)$ for the unique coalgebra $0: 0 \to F0$. The initial algebra and final coalgebra represent finite and infinite traces respectively \cite{rutten2000universal}. For a coalgebra $c: X \to FX$, the algebra $\mu(c)$, may be interpreted as a semantic object for $c$ which represents neither finite or infinite traces. When $F$ is a $\Set$-functor, $\mu(c)$ contains elements of both finite traces and infinite traces. In the colimit for $\mu$, the constants of $F$ generate a copy of its inital algebra. For each element $x \in X$, there is an object of $\mu(c)$ representing its \textit{orbit}. In the colimit for $\mu$ every object is identified with its successor, so in the algebra $\mu(c)$ there is one element for each equivalence class generated by the transitive closure of the sucessor relation. Note that the infinite trace semantics is given by the unique map $c \to \nu(1)$. Transferring this map accross the adjunction gives the unique morphism $\mu(c) \to 1$. To us, this suggests that the algebra $\mu(c)$ is somehow \textit{precompiling} the final trace semantics of $c$. Before moving on to the next section we state a corollary about recursive coalgebras and corecursive algebras. 
\begin{defn}
An $F$-algebra $a$ is corecursive if for every $F$-coalgebra $b$ there is a unique coalgebra to algebra morphism from $b \to b$. Dually, an $F$-coalgebra $b$ is recursive if for any $F$-algebra $a$, there is a unique coalgebra to algebra morphism $b \to a$.
\end{defn} Recursivity of a coalgebra relates to the termination of that coalgebra when thought of as a program (c.f. \cite{adamek_lucke_milius_2007}). Recursivity is closely related to our adjunction.
\begin{cor}
The initial algebra for $F$ is recursive and the final coalgebra for $F$ is corecursive.
\end{cor}
\begin{proof}
 The proof of the adjunction implies that $\Coalg(F)(c,\nu(1)) \cong \mathsf{CoAlgToAlg}(c,1)$ where the latter set is the set of coalgebra to algebra homorphisms into the terminal algebra. This set has a unique element implying that $\nu(1)$ is corecursive. A similar proof holds for the dual statement.
\end{proof}
\section{$\mu$ and $\nu$ Coincide With a Dagger}
When coalgebras for a polynomial functor $F:\Set \to \Set$ are interpreted as $F$-shaped automata, the initial $F$-algebra serves as finite trace semantics and the terminal $F$-coalgebra gives an infinite trace semantics. When $F$ is no longer a $\Set$-functor this interpretation breaks down. For example if $F: \Rel\to \Rel$, where $\Rel$ is the category of sets and relations, then the initial algebra and terminal coalgebra coincide \cite{smyth1982category}. In \cite{karvonen2019way}, it is shown that this holds more generally in any dagger category. With this coincidence, the initial algebra/final coalgebra gives a finite trace semantics instead of an infinite trace semantics. To obtain a semantics for infinite traces, Urabe and Hasuo construct an object which is weakly terminal among coalgebras and define the infinite trace semantics as the maximal map into this object \cite{hasuo2018coalgebraic}. Note that the limit colimit coincidence causes no issues when $\mu(c)$ is interpreted as a semantic object for $c$. However, A generalized limit colimit coincidence also holds for the fixed points generated by $\mu$ and $\nu$.
\begin{defn}
A dagger category $(C,\dag)$ is a category equipped with an identity on objects functor $\dag: C \to C^{\op}$ such that $\dag^2=id$.\end{defn}
\begin{theorem}
    Suppose that $(C,\dag)$ is a dagger category with limits and colimits of countable chains and $F: C \to C$ is a dagger functor preserving such limits and colimits. Then  there is an isomorphism
    \[ \mu(c)^\dag \cong \nu(c^\dag)\]
    for each coalgebra $c$. Dually, for each algebra $a$, there is an isomorphism $\nu(a)^\dag \cong\mu(a^\dag)$.
\end{theorem}
\begin{proof}
For a coalgebra $X \xrightarrow{c} FX$ we have
\begin{align*}
    \nu(c^{\dag}) & \cong \lim(X \xleftarrow{c^{\dag}} FX \xleftarrow{Fc^{\dag}} F^2X \leftarrow \ldots ) \\
    &\cong \colim_{C^{\op}} (X \xleftarrow{c^{\dag}} FX \xleftarrow{Fc^{\dag}} F^2X \leftarrow \ldots ) \\
    &\cong \colim(X \xrightarrow{c} FX \xrightarrow{Fc} F^2X \to \ldots )^{\dag} \\
    &\cong \mu(c)^{\dag}
\end{align*}
The second isomorphism is because limits in $C$ are colimits in $C^{op}$ and the third isomorphism is because $\dag$ preserves colimits because it is an equivalence. A similar proof holds for the dual statement.
\end{proof}
\section{Conclusion}In this extended abstract, we have argued that \textit{middle fixpoints} i.e.\ fixpoints which are neither initial or terminal are interesting enough to merit further study. The adjunction $\mu \vdash \nu$, is closely related to the concept of coalgebra to algebra homomorphisms. These have been studied in the case when they are unique through the notions or recursive algebras and corecursive coalgebras \cite{capretta2009corecursive}. However, in \cite{hauhs2015scientific}, the author argued that coalgebra to algebra morphisms also hold interest when they are not unique. Using examples in probability, dynamical systems, and a game theory, the authors showed how non-unique coalgebra to algebra morphisms often represent solutions to problems in these disciplines. A morphism out of $\mu(c)$ represents a coalgebra to algebra homomorphism originating in $c$ and dually for $\nu(a)$. In future work, we hope to recast the examples given in \cite{hauhs2015scientific} in terms of our adjunction $\mu$ and $\nu$ to further explore properties of these fixpoints as semantics for coalgebras and algebras.

\newcommand{\etalchar}[1]{$^{#1}$}

\end{document}